\documentclass[12pt]{article}

 \usepackage{amsfonts,amssymb,eucal}

\usepackage{amsthm}
 \usepackage{amsmath}
\usepackage{amscd}
 \usepackage{latexsym} 
\numberwithin{equation}{section}
\newcommand{\ad}{{\mbox{\rm ad}}}
 \newcommand{\e}{{\mbox{\rm e}}}

 \newcommand{\mb}[1]{{\mbox{\boldmath{$#1$}}}}
  \newcommand{\mc}[1]{{\mathcal{#1}}}
 \newcommand{\got}[1]{{\mathfrak{#1}}}
\newcommand{\db}[1]{{\mathbb{#1}}}
 
\newcommand{\pa}{\partial}

\newcommand{\C}{\ensuremath{\mathbb{C}}}
\newcommand{\N}{\ensuremath{\mathbb{N}}}

\newcommand{\T}{\ensuremath{S(U(1)\times U(1) \times U(1))}}
 \newcommand{\Hi}{\ensuremath{\mathcal{H}}}
 
 \newcommand{\Hinf}{\ensuremath{\mathcal{H}^{\infty}}}
  \newcommand{\U}{\ensuremath{\mathcal{U}}}
  \newcommand{\g}{\ensuremath{\got{g}}}
\newcommand{\m}{\ensuremath{\got{m}}}
\newcommand{\bb}{\ensuremath{\got{b}}}
\newcommand{\gc}{\ensuremath{\got{g}_{\C}}}
\newcommand{\Ugc}{\ensuremath{\mathcal{U}({\gc}})}
\newcommand{\Ad}{\ensuremath{{\mbox{\rm{Ad}}}}}
\renewcommand{\P}{\ensuremath{\mathbb{P}}}
 \newcommand{\Ph}{\ensuremath{\P (\Hi )}}
\newcommand{\Phinf}{\ensuremath{\P (\Hinf )}}
 \newcommand{\Gras}{\mbox{$G_n({\C}^{m+n})$}}
\newcommand{\DM}{\ensuremath{{\got{D}}_M }}
\newcommand{\AM}{\ensuremath{{\got{A}}_M }}
 \newcommand{\D}{\ensuremath{{\got{D}}}}
\newcommand{\A}{\ensuremath{{\got{A}}}}
\newcommand{\AAA}{\ensuremath{{\db{A}}_M}}
\newcommand{\am}{\ensuremath{{\bf{A}}_M}}%
 
\newcommand{\FSB}{symmetric Fock space }
 \newtheorem{Theorem}{Theorem}
 \newtheorem{Remark}{Remark}

\newcommand{\fl}{\ensuremath{{\mathcal{F}}_{\Hi}}}
\newcommand{\ep}{\ensuremath{{\epsilon}}}
\newcommand{\NC}{\ensuremath{\mathcal{G}}}
\newcommand{\HR}{\ensuremath{\it{H}}}
\newtheorem{Proposition}{Proposition}
 \newtheorem{lemma}{Lemma}
\theoremstyle{definition}

\newtheorem{deff}{Definition}

\begin{document}

\begin{center}
{\Large {\bf Linear Hamiltonians on homogeneous K\"ahler manifolds of
 coherent states}}\\[2ex]
S. Berceanu, A. Gheorghe\\[2ex]
 National 
 Institute for Physics and Nuclear Engineering\\
         Department of Theoretical Physics\\
         PO BOX MG-6, Bucharest-Magurele, Romania\\
         E-mail: Berceanu@theor1.theory.nipne.ro; 
Cezar@theor1.theory.nipne.ro\\
\end{center}

\begin{abstract}

 Representations of coherent state Lie algebras on coherent state
manifolds as first order differential operators are presented.
 The explicit expressions of the differential
action of the generators  of  semisimple Lie groups determine for
 linear Hamiltonians in the generators of the groups first
order differential equations of motion with holomorphic
polynomials coefficients.
 For hermitian symmetric manifolds the equations of motion
are matrix Riccati equations.
It is presented the simplest example of the non-symmetric space
$SU(3)/S(U(1)\times U(1)\times U(1))$ where the polynomials describing
the equations of motion have 
the maximum degree  3.
\end{abstract}
\section{Introduction}
In references \cite{sbcag,sbl} it was shown   that
 the differential action of the
 generators of the groups  on  coherent state manifolds
  which have the structure of
  hermitian symmetric spaces can be written down as a sum
 of two terms, one a polynomial $P$, and the second one  a sum of partial
 derivatives times some  polynomials $Q$-s,
    the degree of  polynomials  being  less than 3.
   Our investigations on the differential action of the
generators has been extended from
 hermitian  groups acting  on hermitian symmetric spaces
to semisimple Lie groups acting   on coherent state manifolds
which admit a K\"ahler structure, and  
 explicit formulas for the polynomials $P$ and $Q$-s  have been given
 \cite{sbcpol}.  Similar investigation has been done in 
 \cite{dob}.
  Explicit formulas for the
simplest example of a compact non-symmetric coherent state     
 manifold, $SU(3)/\T$, where the degree of the
polynomial is already 3, have been also obtained  \cite{sbcpol}. We
have 
formulated the problem of the
 differential action of the generators of the
 so called
coherent state (shortly, CS) groups \cite{lis1,lis2,lis,neeb} in
\cite{sbc2002}.
These are
 groups whose quotient with the stationary groups are manifolds
which admit a holomorphic embedding in a projective Hilbert
space. This class of groups contains all compact groups, all simple
hermitian groups, certain solvable groups and also mixed groups as the
semidirect product of the Heisenberg group and the symplectic group
\cite{neeb}.   

The 
 coherent states are a  useful tool of investigation of quantum and classical
 systems \cite{perG}. It was shown in \cite{sbcag,sbl} that a linear 
Hamiltonian in the
 generators of the groups implies  equivalent quantum and classical
evolution. It  was proved that for Hermitian symmetric spaces the evolution
 equation generated by Hamiltonians which are linear in the generators of the
 group is a matrix Riccati equation.
 It is interesting to see how it looks like
 the corresponding equation of
 motion generated by linear Hamiltonians for CS-manifolds.  This
question is the main topic of the present investigation. 

 Another field of possible applications is the determination of the Berry
   phase \cite{berry} on CS-manifolds.
 In \cite{sbcag,sbl} there were presented explicit expressions for the 
 Berry phase for the complex Grassmann manifold. These results 
were used further 
 \cite{sb7} for explicit calculation of the symplectic area of geodesic
 triangles on the complex Grassmann manifold.

 The paper is laid out as follows. In the first part we recall some
facts
about realization of coherent state Lie algebras by differential
operators
(cf. \cite{sbc2002}).
\S \ref{CSR} contains
the definition of  CS-orbits, 
 in the context of Lisiecki \cite{lis1,lis2,lis} 
and Neeb \cite{neeb}. The geometry of coherent state manifolds for
compact groups was previously considered in \cite{morse}.
\S \ref{CSvectors} deals with the so called
Perelomov's CS-vectors.  In \S \ref{fock} we construct the space of
functions on which the differential operators will act.
In \S \ref{DIFF} we study the representations of Lie algebras of
CS-groups by differential operators.  \S \ref{sSIMPLE} deals with the
semisimple case. \S \ref{liealg} recalls some standard facts about the
semisimple Lie algebras.  Perelomov's coherent state vectors for
semisimple Lie groups are defined in  \S \ref{perv}.  \S \ref{MAIN}
recalls the results established in \cite{sbcpol}. The example
of $SU(3)/\T$ is presented in \S \ref{noul}.  In the second part
the equations of motion associated to linear Hamiltonians in the
generators of the groups  are investigated. Some known facts
established in references \cite{sbcag,sbl} are recalled in \S
\ref{eqm}. The next sections present firstly
 the example of
equations of motion generated by linear Hamiltonians in the generators
of the oscillator 
group (\S \ref{osc}) and on the Riemann sphere and its non-compact dual
(\S\ref{sfera}). The case
 of the complex Grassmann manifold is summarized in \S\ref{grasu}. The
last example in \S\ref{ultim}
 presents the equations of motions generated by the
differential operators of \S \ref{noul}.

We use for the scalar product the convention:
$(\lambda x,y)=\bar{\lambda}(x,y)$, $ x, y\in\Hi ,\lambda\in\C $.

\section{Coherent state representations}\label{CSR}

Let us consider the triplet $(G, T, \Hi )$, where $T$ is
 a continuous, unitary
representation 
 of the  Lie group $G$
 on the   separable  complex  Hilbert space \Hi .
Let us denote by $\Hinf$ the dense subspace of \Hi~ consisting of
 those vectors
$v$ for which the orbit map $G\rightarrow \Hi , g\mapsto
 T (g).v$ is smooth.  Let us
pick up $e_0\in \Hinf$ and let  the notation
$e_{g,0}:=T(g).e_0, g\in G$.
We have an action $G\times \Hinf\rightarrow\Hinf$, $g.e_0 :=
e_{g,0}$. When there is no possibility of confusion, we write just
$e_{g}$ for $e_{g,0}$. 
Let us denote by  $[~]:\Hi^*:=\Hi\setminus\{0\}\rightarrow\Ph=\Hi^*/ \sim
$ the projection with respect to the equivalence relation
 $[\lambda x]\sim [x],~ \lambda\in \C^*,~x\in\Hi^*$. So,
$[.]:\Hi^*\rightarrow \Ph , ~[v]=\C v$. The action 
$ G\times \Hinf\rightarrow\Hinf$ extends to the action
$G\times \Phinf\rightarrow\Phinf , g.[v]:=[g.v]$. 

Let us now denote by $H$  the isotropy group $H:=G_{[e_0]}:=
\{g\in G|g.e_0\in\C e_0\}$.
We shall consider generalized coherent 
 states on complex  homogeneous manifolds $M\cong
G/H$, imposing the restriction that $M$ be  a complex submanifold of
\Phinf .

\begin{deff}\label{def1}
 a) The orbit $M$ is called a CS-{\em orbit} if 
there exits  a holomorphic embedding
$\iota : M \hookrightarrow \Phinf$.
In such a case $M$ is also called CS-{\it manifold}.

b) $(T,\Hi )$ is called a  CS-{\em representation} if there exists a cyclic
 vector $0\neq
e_0\in\Hinf $ such that $M$ is a CS-orbit.

c)  The groups $G$ which admit
CS-representations  are called CS-{\it groups}, and their Lie algebras
\g~ are called CS-{\it Lie algebras}.
\end{deff}

The  $G$-invariant complex structures on the homogeneous
 space $M=G/H$   
 can be introduced in an algebraic manner. 
 For $X\in\g$, where \g ~is the Lie algebra of the Lie group $G$,  let
us define the unbounded operator $dT(X)$ on \Hi~ by
$dT(X).v :=\left. {d}/{dt}\right|_{t=0} T(\exp tX).v$
whenever the limit on the right hand side exists. We obtain a
representation of the Lie algebra \g~ on \Hinf , {\it the derived
representation}, and we denote
${\mb{X}}.v:=dT(X).v$ for $X\in\g ,v\in \Hinf$. Extending $dT$ by complex
linearity, we get a representation of the complex Lie algebra \gc~ on
the complex vector space \Hinf . 
Lemma XV.2.3 p. 651 in \cite{neeb}
and Prop. 4.1 in \cite{lis} determine when a smooth vector
generates a complex orbit in \Phinf .
We denote by $B:=<\exp_{G_{\C}} \bb >$ the Lie group
 corresponding to the Lie algebra
\bb ,  with
$\got{b}:=\overline{\got{b}(e_0)}$, where $\got{b}
 (v):=\{X\in\gc :X.v\in \C v\}=(\g_{\C})_{[v]}$. The group $B$ is
 closed in the complexification $G_{\C}$ of
$G$, cf. Lemma XII.1.2. p. 495 in \cite{neeb}. 
 The  complex structure on $M$ is induced by an 
embedding in a complex manifold,
 $i_1:M\cong G/H \hookrightarrow  G_{\C}/B$.
 We  consider such
manifolds  which admit a holomorphic
 embedding $i_2: G_{\C}/B\hookrightarrow
\Phinf$.   Then the embedding $\iota
=i_1\circ i_2$,  $\iota : M  \hookrightarrow \Phinf$
 is a holomorphic embedding,  and  the complex structure comes as in 
  Theorem XV.1.1 and Proposition XV.1.2 p. 646 in \cite{neeb}.

\section{Coherent state vectors }\label{CSvectors}

Now we construct what we  call  Perelomov's
generalized  coherent state vectors, or simply CS-vectors,  based on the
CS-homogeneous manifolds $M\cong G/H$.

  We denote also  by $T$ the holomorphic extension of the
representation $T$ of $G$ to  the complexification $G_{\C}$ of $G$,
 whenever this
holomorphic extension exists. In fact, it can be shown that
 in the situations under interest
in this paper, this holomorphic extension exists 
\cite{neeb95,neeb96}.
 Then there exists the 
homomorphism $\chi_0 $ ($\chi$), $\chi_0: ~H\rightarrow \db{T}$,
($\chi : B\rightarrow \C^*$), 
such that 
$ H = \{g\in G| e_g=\chi_0(g)e_0\}$ (respectively, 
$B = \{g\in G_{\C}| e_g=\chi(g)e_0\})$, where
$\db{T}$ denotes the torus $\db{T}:=\{ z\in\C| |z|=1\}$. 

For the homogeneous space $M=G/H$ of cosets $\{gH\}$, let $\lambda
:G\rightarrow G/H$ be the natural projection $g\mapsto gH$, and let
$o:=\lambda (\mb{1})$, where $\mb{1}$ is the unit element of $G$. 
 Choosing a section $\sigma :G/H\rightarrow
G$ such that $\sigma ( o )={\bf{1}} $, every element $g\in G$ can
be written down as $g=\tilde{g}(g)h(g)$, where
$\tilde{g}(g)\in G/H$ and $h(g)\in H$. Then we have
$e_g=e^{i\alpha (h(g))}e_{\tilde{g}(g)}$,
where $e^{i\alpha(h(g))}=\chi_{0}(h)$. Now we take into account that
$M$ also admits an embedding in  $G_{\C}/B$. We choose a local system
of coordinates parametrized by $z_g$ (denoted also simply $z$, where
there is no possibility of confusion) on  $G_{\C}/B$. Choosing a
section 
$G_{\C}/B
\rightarrow G_{\C}$ such that any element $g\in G_{\C}$ can be written
as $g=\tilde{g}_bb(g)$, where $\tilde{g}_b\in G_{\C}/B$, and
$ b(g)\in B$, we have
$e_g=\Lambda (g)e_{z_g}$,
where
$\Lambda (g)= \chi
(b(g))=e^{i\alpha(h(g))}(e_{z_g},e_{z_g})^{-\frac{1}{2}}$.

 Let    us denote by
$\got{m}$
 the vector space orthogonal to $\got{h}$ 
 of
the Lie algebra \g , i.e. we have the vector space
decomposition $\g=\got{h}+\got{m}$. Even more, it can be shown that the
 vector space decomposition $\g=\got{h}+\got{m}$ is
Ad $H$-invariant. The homogeneous spaces $M\cong G/H$ with this
decomposition are called 
 {\em reductive spaces} (cf. \cite{nomizu})
and it can be proved that {\em the
CS-manifolds  are  reductive spaces} (cf. \cite{sbc2002}).
So, {\it the tangent space to $M$ at $o$
can be identified with $\got{m}$}.
Recall 
 that for CS-groups the
 CS-representations are highest weight representations
  and the vector
$e_0$ is a primitive element of the generalized  parabolic algebra
$\got{b}$  (cf. \cite{neeb}).

Let us denote ${\mb{X}}:=dT(X), X\in \Ugc$, where    $\U$  denotes the
universal enveloping algebra. Let 
$\tilde{g}(g)=\exp X, \tilde{g}(g)\in G/H,~ X\in\got{m}$,
$e_{\tilde{g}(g)}=\exp({\mb{X}})e_0$.
Let us remember again   
 Theorem XV.1.1 p. 646 in \cite{neeb}.
Note  that $T_o(G/H)\cong\got{g}/\got{h}\cong \gc/\bar{\got{b}}\cong
(\got{b}+\bar{\got{b}})/\bar{\got{b}}\cong
\got{b}/h_{\C}$,
 where we
have a linear isomorphism $\alpha:\got{g}/\got{h}\cong
\gc/\bar{\got{b}}$, $\alpha(X+\got{h})=X+\bar{\got{b}}$ (cf. \cite{neeb1}). 
We can take
instead of $\got{m}\subset\g$ the subspace $\got{m}' \subset \gc$
complementary to $\bar{\bb}$, or the subspace of $\got{b}$
complementary to $\got{h}_{\C}$.
If we choose a local {\em canonical} system of coordinates
$\{z_{\alpha}\}$ with respect to the basis $\{X_{\alpha}\}$ in \m',
then we can introduce the vectors
\begin{equation}\label{cvect}
e_{z}=\exp(\sum_{X_{\alpha}\in\got{m}'}z_{\alpha}{\mb{X}}_{\alpha}).e_0
\in\Hi .
\end{equation}
We get
\begin{equation}\label{unu2}
e_{\sigma (z)}=T(\sigma (z)), ~z\in M, 
\end{equation}
 and  we prefer to choose local coordinates such that
\begin{equation}\label{doi}
e_{\sigma (z)}=N(z)e_{\bar{z}},~~
N(z)=(e_{\bar{z}},e_{\bar{z}})^{-1/2}.
\end{equation}
 Equations (\ref{cvect}), (\ref{unu2}), and (\ref{doi})
 define locally  the {\em coherent vector
 mapping}
\begin{equation}\label{cvm}
\varphi : M\rightarrow \bar{\Hi}, ~ \varphi(z)=e_{\bar{z}},  
\end{equation}
where $ \bar{\Hi}$ denotes the Hilbert space conjugate to $\Hi$.
We call the  vectors $e_{\bar{z}}\in\bar{\Hi}$ indexed by the points
 $z \in M $  {\it
Perelomov's coherent state vectors}.

\section{ The \FSB ~  \fl }\label{fock}

We have considered homogeneous CS-manifolds $M\cong G/H$ whose complex
structure comes from the embedding $i_1:M\hookrightarrow
G_{\C}/B$. We have chosen a section $\sigma: G_{\C}/B\to  G_{\C}$,
 and $G_{\C}$ can be regarded
as a complex analytic principal bundle     
$ B\stackrel{i}{\rightarrow}G_{\C}\stackrel{\lambda}{\rightarrow}
G_{\C}/B$.

 Let us introduce the function $f'_{\psi}:G_{\C}\to \C$,
$f'_{\psi}(g):=(e_g,\psi ), g\in G, \psi\in \Hi$.
Then 
$f'_{\psi}(gb)=\chi(b)^{-1}f'_{\psi}(g), g\in G_{\C}, b\in B$,
where  $\chi$ is  the continuous homomorphism of the 
isotropy subgroup $B$ of  $G_{\C}$  
  in $\C^*$.   
 {\it The coherent states realize the space
of holomorphic global sections} $
\Gamma^{\text{{hol}}}(M,L_{\chi})=H^0(M,L_{\chi})$ {\it on the}
 $G_{\C}$-{\it homogeneous line bundle
$L_{\chi}$ associated by means of the character
 $\chi$ to the  principal B-bundle}
 (cf. \cite{onofri}).
 The  holomorphic line bundle is $L_{\chi}:=M\times_{\chi}\C$,  
 also denoted   $L:=M\times_B\C$ (cf. \cite{bott,tirao}).

The local trivialization of the line bundle $L_{\chi}$ associates to
every $\psi\in\Hi$ a holomorphic function $f_{\psi}$   on a  open
set  in $M \hookrightarrow G_{\C}/B$. Let the notation $G_S:=G_{\C}
\setminus S$,
 where $S$ is the   set
$S:=\{g\in G_{\C}|\alpha_g=0\}$,
and $\alpha_g:= (e_g,e_0)$. $G_S$ is a dense subset of $G_{\C}$.
We  introduce  the function $f_{\psi}:G_{S}\mapsto\C$,
$f_{\psi}(g)=\frac{f'_{\psi}(g)}{\alpha_g}, \psi\in\Hi ,~g\in G_{S}$.
The
function $f_{\psi}(g)$ on $G_S$
 is actually a function of the projection
$\lambda (g)$,
 holomorphic in  $M_S:=\lambda ( G_S)$.
We have introduced the function 
$f_{\psi}(g)=f_{\psi}(z_g)=\frac{(e_{\bar{z}_{g}},\psi )}
{(e_{\bar{z}_{g}},e_0)}$, where
$(e_{z_{g}},e_0)\not= 0$,
and also the coherent state map
$\varphi :M\rightarrow\overline{\Hi}^{\infty} , \varphi (z)=e_{\bar{z}},
z\in{\mc{V}}_0,$
where the canonical coordinates
$z=(z_1,\ldots ,z_n )$ constitute
 a   local chart   on  ${\mathcal{V}}_0:=M_S\rightarrow
\C^n$,
such that $0=(0,\ldots ,0)$  corresponds to $\{ B\}$.
Note  also that
 ${\mathcal{V}}_0\equiv M\setminus \Sigma_0$, where $\Sigma_0:=\lambda
 (S)$ is the set of points of $M$ for which the
coherent state vectors are orthogonal to $e_0\in \Hi$, called 
{\it
polar divisor} of the point $z=0$  (cf. \cite{sbcl}).

Supposing that
{\it the line bundle $L_{\chi}$ is already very ample},
\fl~ is defined as the set of functions corresponding to sections such
that
$\{f\in L^{2}(M,L)\cap{\mc{O}}(M,L)
|(f,f)_{\fl}<\infty \}$
 with respect to  the scalar product
\begin{equation}\label{scf}
(f,g)_{\fl} =\int_{M}\bar{f}(z)g(z)d\nu_M(z,\bar{z}),
\end{equation} 
where  $d{\nu}_{M}(z,\bar{z})$ is the invariant measure 
$\frac{d{\mu}_M(z,\bar{z})}{(e_{\bar{z}},e_{\bar{z}})}$, and
 $d{\mu}_M(z,\bar{z})$ represents the $G$-invariant
Radon measure on $M$.
It can be shown that the space of functions
 $\fl$ identified with $ L^{2,{\text{hol}}}(M,L_{\chi})$ {\em  is a closed
subspace of $L^2(M,L_{\chi})$ with continuous point evaluation}
(cf. \cite{PWJ}). 
 
Note that eq.  (\ref{scf})
is nothing else than the {\em {Parseval   overcompletness
 identity}}  \cite{berezin}:
\begin{equation}\label{orthogk}
(\psi_1,\psi_2)=\int_{M=G/K}(\psi_1,e_{\bar{z}})(e_{\bar{z}},\psi_2)
d{\nu}_{M}(z,\bar{z}),~ (\psi_1,\psi_2\in \Hi ) .
\end{equation}

Let us  now introduce the map
\begin{equation}\label{aa}
\Phi :\Hi^{\star}\rightarrow \fl ,\Phi(\psi):=f_{\psi},
f_{\psi}(z)=\Phi(\psi )(z)=(\varphi (z),\psi)_{\Hi}=(e_{\bar{z}},\psi)_{\Hi},~
z\in{\mathcal{V}}_0,
\end{equation}
where we have identified the space $\overline{\Hi}$  complex conjugate 
 to \Hi~  with the dual
space
$\Hi^{\star}$ of $\Hi$.

In fact,{\em{ our supposition that $L_{\chi}$ is already a very ample line
bundle implies the validity of eq. (\ref{orthogk})}} (cf. 
 Theorem XII.5.6 p. 542 in \cite{neeb},
 Remark VIII.5 in \cite{neeb94}, and 
 Theorem XII.5.14 p. 552 in \cite{neeb}).

It can be seen 
that  the group-theoretic relation (\ref{orthogk}) on homogeneous
manifolds fits into 
 Rawnsley's global realization \cite{raw} of Berezin's coherent states on 
quantizable K\"ahler manifolds \cite{berezin}.

It can be defined a function  $K$,
   $K: M\times\overline{M}\rightarrow \C$, which on  ${\mathcal{V}}_0\times
\overline{\mathcal{V}}_0$ reads
\begin{equation}\label{kernel}
K(z,\overline{w}):=K_w(z)=
(e_{\bar{z}},e_{\bar{w}})_{\Hi}.
\end{equation}

Taking into account (\ref{aa}) and supposing that
eq. (\ref{orthogk})
is true, it follows that  the function $K$ (\ref{kernel}) is a reproducing
kernel. 
We have:
\begin{Proposition}\label{realization}
Let $(T,\Hi)$ be a CS-representation and let us consider
 the Perelomov's CS-vectors defined in 
(\ref{cvect})-(\ref{doi}). Suppose that the  line bundle $L$ is
 very ample.  Then

\mbox{\rm{i)}} The function $K:M\times\overline{M}\rightarrow \C$, 
 $K(z,\overline{w})$ defined by equation   
(\ref{kernel}) is a   reproducing kernel.

\mbox{\rm{ii)}}  Let \fl~ be
the space $L^{2,\text{\em{hol}}}(M, L)$ endowed with the scalar
 product (\ref{scf}). Then 
 \fl~ is the reproducing kernel Hilbert space $\Hi_K\subset \C^M$
  associated to the kernel $K$
(\ref{kernel}).

\mbox{\rm{iii)}} The evaluation  map $\Phi$ defined in
 eqs. (\ref{aa}) 
extends to an  isometry   
\begin{equation}\label{anti}
(\psi_1,\psi_2)_{\Hi^{\star}}=(\Phi (\psi_1),\Phi
(\psi_2))_{\fl}=(f_{\psi_{1}},f_{\psi_{2}})_{\fl}=
\int_M\overline{f}_{\psi_1} (z)f_{\psi_2}(z)d\nu_M(z),
\end{equation}
and  the overcompletness eq. (\ref{orthogk}) is verified. 
\end{Proposition}

\section{Representations of coherent state Lie algebras  by
 differential operators}\label{DIFF}

We remember the definitions of the functions $f'_{\psi}$  and
$f_{\psi}$,  which allow to
write down
\begin{equation}\label{iar}
f_{\psi}(z)=(e_{\bar{z}},\psi)=\frac{(T(\bar{g})e_0,\psi)}
{(T(\bar{g})e_0,e_0)},~ z\in M,~ \psi\in \Hi .
\end{equation}
We get 
\begin{equation}\label{iar1}
 f_{T(\overline{g'}).\psi}(z)= \mu
(g',z)f_{\psi}(\overline{g'}^{-1}.z),
\end{equation}
where 
\begin{equation}\label{iar2}
 \mu (g',z)=
\frac{(T(\overline{g'}^{-1}\overline{g})e_0,e_0)}{(T(\overline{g})e_0,e_0)}
=\frac{\Lambda (g'^{-1}g)}{\Lambda (g)}.
\end{equation}
Recall that 
${T(g).e_0=e^{i\alpha (g)}e_{\tilde{g}}=\Lambda (g) e_{z_{g}}}$
where we have used the decompositions
$g=\tilde{g}.h, ~(G=G/H.H);~~ g = z_g. b ~(G_{\C}=G_{\C}/B.B)$.
 We have also  the relation
$\chi_0(h)=e^{i\alpha (h)},~ h\in H$ and 
$ \chi (b) =\Lambda (b),~b\in B$, where
$\Lambda(g)=\frac{e^{i\alpha(g)}}{(e_{\bar{z}},e_{\bar{z}})^{1/2}}$.
We can also write down another expression for multiplicative factor
 $\mu$ appearing  in eq. (\ref{iar1})    using the CS-vectors:
\begin{equation}\label{iar3}
\mu (g',z)=\Lambda(\bar{g'})(e_{\bar{z}},e_{\bar{z'}})=e^{i\alpha
  (\bar{g'})}
\frac{(e_{\bar{z}},e_{\bar{z'}})}{(e_{\bar{z'}},e_{\bar{z'}})^{1/2}}.
\end{equation}

The following assertion is easily  checked out:
\begin{Remark}Let us consider the relation (\ref{iar}). Then we have
(\ref{iar1}), where  $\mu$ can be written down as in equations
(\ref{iar2}), 
(\ref{iar3}). 
We have the relation
 $\mu (g,z) =J(g^{-1},z)^{-1}$, i.e. the multiplier $\mu$
 is the cocycle in the  unitary representation
$(T_K,\Hi_K)$ attached to the
positive definite holomorphic kernel $K$ defined by equation (\ref{kernel}),
 \begin{equation}\label{num}
(T_K(g).f)(x):=J(g^{-1},x)^{-1}.f(g^{-1}.x),
\end{equation}
and the cocycle verifies the relation
\begin{equation}\label{prod}
J(g_1g_2,z)=J(g_1,g_2z)J(g_2,z).
\end{equation}
\end{Remark}
  Note that
{\em   the prescription (\ref{num})
 defines
a continuous action of $G$ on ${\mbox{\rm{Hol}}}(M,\C )$ with respect to
the compact open topology on the space  ${\mbox{\rm{Hol}}}(M,\C )$.
 If $K:M\times M\rightarrow \C$ is a continuous positive definite kernel
holomorphic in the first argument satisfying
$K(g.x,\overline{g.y})=J(g,x)K(x,\overline{y})J(g,y)^*$,
$g\in G$, $x,y\in M$, then the action of $G$ leaves the reproducing kernel
Hilbert space $\Hi_K\subseteq {\mbox{\rm{Hol}}}(M,\C )$ invariant and defines
a continuous unitary representation $(T_K,\Hi_K)$ on this space}
(cf. Prop. IV.1.9 p. 104 in   \cite{neeb}).

Let us consider the triplet $(G, T, \Hi )$. Let
 $\Hi^0:=\Hinf$, considered as a pre-Hilbert space, and
let  $B_0(\Hi^0)\subset \mc{L}(\Hi )$ denote the set
  of linear operators
 $A:\Hi^0\rightarrow \Hi^0$
which have  a formal adjoint  $A^{\sharp}:\Hi^0\rightarrow\Hi^0$, i.e.   
$(x,Ay)=(A^{\sharp}x,y)$ for all $x,y\in \Hi^0$.
Note that if  $B_0(\Hi^0)$ is the set of unbounded operators on \Hi ,
then the domain $\mathcal{D}(A^*)$ contains $\Hi^0$ and $A^*\Hi^0\subseteq 
 \Hi^0$, and  it makes sense to refer to the closure $\overline{A}$
 of $A\in B_0(\Hi^0)$ (cf. \cite{neeb} p. 29; here $A^*$ is the
 adjoint of $A$).

Let \g~ be the Lie algebra of $G$ and let us denote
 by $\mathcal{S}:=\Ugc$ the semigroup associated with
 the universal enveloping  algebra equipped with the antilinear involution
 extending the antiautomorphism $X\mapsto X^*:=-\bar{X}$ of $\gc$.  {\it
 The derived representation} is defined as
\begin{equation}\label{derived}
dT:\Ugc\rightarrow B_0(\Hi^0),~~ \mbox{\rm{with}} ~~ dT(X).v
:=\left.\frac{d}{dt}\right|_{t=0}T(\exp tX).v, X\in\g .
\end{equation}
Then $dT$ is a {\it hermitian representation}
 of $\mathcal{S}$ on $\Hi^0$ (cf. Neeb \cite{neeb}, p. 30).
Let us denote his image in $B_0(\Hi^0)$ with $\am := dT(\mathcal{S})$. 
If $\Phi : \Hi^{\star}\rightarrow \fl $ is the  isometry
(\ref{aa}), we are interested in the study of the image of 
\am~  via $\Phi$  as subset in the
algebra of holomorphic, linear differential operators,  
$ \Phi\am\Phi^{-1}:={\db{A}}_M\subset\got{D}_M$.

The {\it  sheaf}
 $\got{D}_M$ (or simply \D ) {\it of holomorphic, finite order, linear
differential operators on} $M$ is a
 subalgebra of homomorphism ${\mathcal Hom}_{\C}({\cal O}_M,{\cal O}_M)$ 
 generated
 by the sheaf ${\cal O}_M$ of germs of holomorphic functions of $M$ and the
 vector fields. 
 We consider also {\it  the subalgebra} \AM~ of ${\db{A}}_M$~
 {\it of differential operators with
 holomorphic polynomial coefficients}.
Let $U:=\mathcal{V}_0$ in $M$, endowed with the 
coordinates
$(z_1,z_2,\cdots ,z_n)$. We set $\pa_i:=\frac{\pa}{\pa z_i}$ and
$\pa^{\alpha}:=\pa^{\alpha_1}_1 
\pa^{\alpha_2}_2\cdots \pa^{\alpha_n}_n$, $\alpha :=(\alpha_1,
\alpha_2 ,\cdots ,\alpha_n)\in\N^n$. The sections of \DM~ on $U$ are
$A:f\mapsto \sum_{\alpha}a_{\alpha}\pa^{\alpha}f$,
$a_{\alpha}\in\Gamma (U, {\cal{O}})$, the $a_{\alpha}$-s being zero
except a finite number. 

For $k\in\N$, let us denote by $\D_k$ the subsheaf of differential
operators of degree $\le k$ and by $\D'_k$ the subsheaf of elements of 
$\D_k$ without constant terms.  $\D_0$ is identified with
$\cal{O}$ and $\D'_1$ with the sheaf of vector fields. The
filtration of \DM~ induces a filtration on  ${\A}_M$.

Summarizing, we have the following three objects which corresponds
each to other:
\begin{equation}\label{correspond}
\g \ni X \mapsto\mb{X}\in\am\mapsto\db{X}\in\AAA\subset \DM, {\text
  {~differential  operator  on}}~ \fl .
\end{equation}
We can see that
 \begin{Proposition}If $\Phi$ is the isometry
(\ref{aa}), then 
$\Phi dT(\g_{\C})\Phi^{-1}\subseteq \D_1$.
\end{Proposition}
{\it Proof}. Let us consider an element in $\g_{\C}$ and his image in \DM 
via the correspondence (\ref{correspond}):
$$\g_{\C}\ni G \mapsto{\db{G}}\in \DM;~~
 {\db{G}}_z(f_{\psi}(z))= 
\db{G}_z(e_{\bar{z}},\psi)= (e_{\bar{z}},\mb{G}\psi),$$
$$\mb{G}=dT (G)=\frac{d}{dt }|_{t=0}T(\exp (tG)).$$
Remembering equation
(\ref{iar1}) and  determining  the derived representation,
 we get
\begin{equation}\label{sss}
\db{G}_{z}(f_{\psi}(z))=(P_{G}(z)+
\sum Q^i_G(z)\frac{\pa}{\pa  z_i})f_{\psi}(z);~
\end{equation}
 $$P_G(z)=\frac{d}{dt}|_{t=0}\mu (\exp (tG),z);~ 
 Q^i_G(z)=\frac{d}{dt}|_{t=o}(\exp(-tG).z)_i . $$
Now we formulate the following assertion:
\begin{Remark}\label{main}
If $(G,T)$ is a CS-representation, then \AAA~ is a subalgebra of
holomorphic differential operators with polynomial coefficients,
 $\AAA \subset \AM\subset \DM$. 
More exactly, for    $X\in\g$,  
  let us denote by $\mb{X}:= dT(X)\in\am $, where the action is
considered on the space of functions \fl . Then, for 
CS-representations, $\db{X}\in \A_1=\A_0\oplus \A_1'$.
 
Explicitly, if  $\lambda\in\Delta $ is a root and
 $G_{\lambda}$ is in a base of the Lie algebra $\g_{\C}$ of $G_{\C}$, then
 his image  ${\db G}_{\lambda}\in\DM$ acts as a first order
 differential operator on the
\FSB \fl
  \begin{equation}\label{generic}
{\db G}_{\lambda}= P_{\lambda}+\sum_{\beta \in\Delta_{\m '}} Q_{\lambda
 ,\beta} \pa_{\beta},~ \lambda \in \Delta,
 \end{equation}
where $ P_{\lambda}$ and $ Q_{\lambda ,\beta}$  are 
 polynomials in
 $z$, and $\m '$ is the subset of $\g_{\C}$ which appears in the
definition  (\ref{cvect}) of the coherent state vectors.
\end{Remark}
 Actually, we don't have a proof of this 
assertion for the general case of CS-groups.
  For the compact case, there exists the calculation of
Dobaczewski \cite{dob}. For compact hermitian symmetric
spaces it was shown \cite{sbcag} that degrees of the polynomials
 $P$ and $Q$-s are   
 $\le 2$ and similarly for the non-compact hermitian symmetric  case
\cite{sbl}.
Neeb \cite{neeb} gives a proof of this Remark
  for 
CS-representations for  the (unimodular) Harish-Chandra type
groups.
{\it If $G$ is an admissible Lie group such that the universal
complexification $G\rightarrow G_{\C}$ is injective and $ G_{\C}$ is
simply connected, then $G$ is of Harish-Chandra type}
(cf. Proposition V.3 in \cite{neeb94}). The derived
 representation (\ref{derived}) is obtained
   differentiating eq. (\ref{num}), and we get two 
 terms, one in $\D_0$ and the
other one in $\D_1'$.
A proof that the two parts are in fact $\A_0$ and respectively $\A_1'$
    is contained  in
Prop. XII.2.1 p. 515 in \cite{neeb}
  for the groups of Harish-Chandra type in the particular situation
where the space $\got{p}^+$ in  Lemma VII.2.16 p. 241 in \cite{neeb}
 is abelian.

\section{Representation of  semisimple  Lie groups by differential
  operators}\label{sSIMPLE} 


\subsection{ Semisimple  Lie groups and flag manifolds}\label{liealg}

We use standard notation referring to Lie algebras of a complex
 semisimple Lie group $G$  
 \cite{wolf}. In this case $\Delta\equiv\Delta_s$, i.e. $\Delta_r=\{\emptyset\}$,
 i.e. all roots are semisimple.

 ${\got g}$ -- complex semisimple Lie algebra

 ${\got t}\subset\got{g}$ -- Cartan subalgebra

 ${\got b}={\got t}+{\got b}^{u}$ -- Borel subalgebra

 ${\got b}^{u}=\sum_{\alpha\in\Sigma^+}{\got g}_{\alpha}$ -- the nilradical of
 ${\got b}$

 $\Sigma$ --  root system for $({\got g} ,{\got t})$
  
$\Sigma^+$ -- a positive root system
 
$\Psi $ -- a simple root system for  $\Sigma$

 $\Sigma\ni\alpha = \sum _{\mu\in\Psi}n_{\mu}(\alpha )\mu$ -- unique,
 $ n_{\mu}\in \N,
 n_{\mu} (\alpha )\ge 0~{\mbox{\rm {if~}}} \alpha \in \Sigma^+;
 n_{\mu} (\alpha )\le 0~{\mbox{\rm {if~}}}  \alpha \in \Sigma^-$ 

  $ \Psi \supset \Phi \longrightarrow
 \Phi^r =\{\alpha \in\Sigma ; n_{\mu}(\alpha )=0$
 whenever $\mu \notin \Phi\}$

 $ \Phi^u =\{\alpha \in\Sigma ; n_{\mu}(\alpha )>0$
 for some  $\mu \notin \Phi\}=\Sigma^+\setminus\{\Sigma^+\cap\Phi^r\}$

 ${\got p}_{\Phi}={\got p}^r_{\Phi}+
{\got p}^{u}_{\Phi}$- parabolic subalgebras of ${\got
 g}$ corresponding to $\Psi\subset\Phi$

 ${\got p}^r_{\Phi}={\got t}+\sum_{\alpha\in\Phi^r}{\got g}_{\alpha}$ -- 
the reductive part of ${\got p}_{\Phi}$

  ${\got p}^{u}_{\Phi}=\sum_{\alpha\in\Phi^u}{\got g}_{\alpha}$
-- the unipotent part of ${\got p}_{\Phi}$ 

 $\Delta_0=\Phi^r; ~\Delta_-=-\Phi^u; ~\Delta_+=\Phi^u$

 $B=\{g\in G; \Ad (g){\got b}={\got b}\}$ -- Borel subgroup
 (maximal solvable)

 $P=\{g\in G; \Ad(g){\got p}={\got p}\}$ --  parabolic subgroup (contains  a
 Borel subgroup).

In the notation of Definitions
 VII.2.4 p. 234, VII.2.6 p. 236 and VII.2.22, p. 244 in \cite{neeb}
 we have $\Phi^u \equiv\Delta^+_p$
 and
$\Phi^r\equiv\Delta_k$.

 We remember the following facts:

1. {\it Let $G$ be complex semisimple Lie group, $G_u$ his compact real
form. Then the isotropy group  $K= G_u\cap P$ 
is connected and the compact simply connected K\"ahler manifold 
$M\approx G_u/K\approx G/P$ is an algebraic
manifold (Hodge),} called generalized complex flag manifold. $G$ is viewed
as a group of holomorphic transformation on $M$.   
 
2.  Let
 {\it $ G $ be a homogeneous compact K\"ahler manifold. Then
 the projective space orbit of an extreme weight vector
is an
 irreducible finite representation of $ G^{\C}$}.

3. {\em 
 Let $G_0$ be a real form of $G, x\in G/P$. Then $K=G_0\cap P=G_u\cap P$. 
  $G_0(x)$ has a
 $G_0$-invariant K\"ahler metric. The K\"ahler orbit $G_0(x)$ is open in
 $G^{\C}/P$}. 

For symmetric spaces $K$ is a maximal compact subgroup of $G$. 
  Hermitian symmetric spaces correspond to centerless
semisimple Lie groups, which verify
the condition  $\got{z}_{\got{g}}(\got{z}(\got{k}))$ $=\got{k}$ of Definition
 VII.2.15 p. 241 in \cite{neeb}) of quasihermitian groups,
and $\got{p}^u_{\Phi}$ is abelian.

We  also need  the commutation relations in the Cartan-Weyl basis \cite{helg}
\begin{equation}\label{cw}
 \left\{
 \begin{array}{ccll}
 \left[ H_i,H_j \right]   & = &   0, & i=1,\dots ,r, H_i\in\got{t},\\
  \left[ H_i,E_{\alpha}\right]  & = &\alpha_i E_{\alpha}, &  \alpha_i =
\alpha (H_i) ,\\
  \left[ E_{\alpha},E_{\beta}\right]  & =& 
 n_{\alpha ,\beta}E_{\alpha +\beta}, & \alpha +\beta \in\Delta\setminus
 \{0\}  ,\\
\left[ E_{\alpha},E_{\beta}\right]  & =& 
 0, & \alpha +\beta \notin\Delta\cup
 \{0\},\\
 \left[ E_{\alpha} ,E_{-\alpha}\right]  & =&
 H_{\alpha}=\sum\alpha_iH_i. & 
\end{array} 
\right.
 \end{equation}
As a consequence, we have also  the commutation relations:
  \begin{equation}\label{cw1}
 \left\{
 \begin{array}{l}
 \left[E_{-\gamma},E_{\gamma}\right]=-\gamma H,~
\gamma H : =  (\gamma , H)  = \sum ^r_{j=1}\gamma_jH_j; \\
 \left[H,E_{\alpha}\right]=\alpha (H) E_{\alpha}.
\end{array}
 \right.
\end{equation}

\subsection{ Perelomov's coherent  vectors for semisimple 
Lie groups}\label{perv}

All representations of compact Lie groups are CS-representations
because these representations are highest weight representations. 
Kostant and Sternberg \cite{ks} showed that for any representation of
a compact group $G$ the orbit to a projectivized highest weight vector is the 
only K\"ahler coherent state orbit.
Harish-Chandra \cite{hc} has defined highest weight representations
for non-compact semisimple (or even reductive) Lie groups. He has classified
 square integrable highest weight representations. This classification has
been fully realized by Enright, Howe and Wallach, and independently
by  Jakobsen \cite{eww}.
  Lisiecki has emphasized  (cf. \cite{lis1} and
  Theorem 6.1 in \cite{lis}) that:
 {\em a non-compact semisimple Lie group is a CS-group if and 
only if it is hermitian. If this is the case, the
$CS$-representations of $G$ are precisely the highest weight representations. 
Each of them has a unique $CS$-orbit, which is the orbit trough the highest
 line}. The starting point of the proof of
Lisiecki is the paper of Borel \cite{bor}, where it is proved: {\em
a non-compact semisimple Lie group $G$ admits a homogeneous
K\"ahler orbit if and only if it is  hermitian, and such a
manifold is of the form $G/Z_{G(S)}$, where $Z_{G(S)}$ is the centralizer of
a torus $S\subset G$; moreover,    it is a holomorphic fiber bundle over the
Hermitian symmetric space $G/K$, where $K$ is a maximal compact
subgroup of $G$,
 with (compact) flag manifolds $K/Z_{G(S)}$
as fibers}.

Let us consider again the triplet $(G,\pi,\Hi )$ where $(G,\pi)$ is a
CS-representation. Then this representation can be realized as an
extreme weight representation. For linear connected 
reductive groups with $Z_K(\got{z})=\got{k}$, where
 $\got{z}$ denotes the center of the
Lie algebra $\got{k}$ of $K$, the effective representation is
furnished by the Harish-Chandra theorem 
(cf. e.g. \cite{knapp}, p. 158).  The theorem furnishes the holomorphic
discrete series for the non-compact case, and  for the compact case it
is 
equivalent with the Borel-Weil theorem  (\cite{bw}; also
cf. \cite{knapp}, p. 143).

In accord with the procedure of \S \ref{CSvectors} for getting Perelomov's
  CS- vectors, we start with $e_{g,0}$. 
Then
  we consider for a everywhere dense subset $G^0\subset G$ the
Gauss decomposition 
\begin{equation}\label{gauss}
g= g_+.g_0.g_- ~,g\in G^0.
\end{equation}
If we restrict ourself to the complex semisimple Lie groups, then, in
the notation of the  \S \ref{liealg}, eq. (\ref{gauss}) reads
\begin{equation}\label{gauss1}
g=\exp\sum_{\alpha\in\Delta_+}\bar{z}_{\alpha}E_{\alpha}.
 \exp\sum_{i=1}^{r}{v}_{i}H_{i}.
\exp\sum_{\alpha\in\Delta_0}{\xi}_{\alpha}E_{\alpha}.
\exp\sum_{\alpha\in\Delta_-}{y}_{\alpha}E_{\alpha},
\end{equation}
where ($\bar{z}_{\alpha},v_i,\xi_{\alpha}, y_{\alpha}$) are local
coordinates for $G$.

If the extreme weight $j$ (here minimal) of the representation has the 
components
 $j = (j_1,\cdots ,j_r)$, 
  where $ r$ is the  rank 
 of the Cartan algebra,  then
 \begin{equation}\label{acth}
 \left\{
 \begin{array}{l}
 \mb{H}_k \e_{j}  =  j_k \e_{j}, k=1,\dots , r ;\\
 \mb{E}_{\alpha}\e_{j}  =  0: \alpha \in \Delta_-\cup  \Delta_0.
\end{array}
 \right.
\end{equation}
 Perelomov's CS-vectors are
\begin{equation}\label{t}
 e_{z,j}=\exp (\sum_{\alpha\in\Delta_+}z_{\alpha}{\mb{E}}_{\alpha  })e_j,
  \end{equation}
where $z_{\alpha}$ are local coordinates for the coordinate
 neighborhood
  ${\mathcal{V}}_0\subset M$, 
  ${\mathcal{V}}_0=M\setminus \Sigma_0$. $\Sigma_0$ is the polar divisor
$\Sigma_0=\lambda (S)$.
 For hermitian symmetric 
spaces it was proved in  \cite{sbcl} that $\Sigma_0= CL_0$, where $CL_0$ is
 the cut locus relative to the origin $0\in M$.

\subsection{Differential operators on semisimple Lie group orbits}\label{MAIN}

We start  introducing the   notation
 $$ Z=\sum_{\alpha\in\Delta_+}z_{\alpha}\mb{E}_{\alpha}, $$
and then $\partial_{\alpha}(Z)=\mb{E}_{\alpha}$, where  $ \partial_{\alpha}:=
 \frac{\partial}{\partial z_{\alpha}},~ \alpha \in \Delta_+ $.

We remember the definition of the  Bernoulli numbers $B_i$ (see e.g.
\cite{abr}):
 
\begin{equation}
 \label{defc}\frac{x}{1-\e^{-x}}=1 +\frac{1}{2}x+\sum_{k\ge 1}(-1)^{k-1}
 \frac{B_k x^{2k}}
 {(2k)!}=\sum_{n\ge 0}c_nx^n ;
 \end{equation}

 \begin{equation}\label{cb}
 c_0=1;~~c_1=\frac{1}{2};
~~c_{2k+1}=0;~~c_{2k}=\frac{(-1)^{k-1}}{(2k)!}B_k ;
 \end{equation}

 \begin{equation}\label{bernuli}
 B_1=\frac{1}{6};~~B_2=\frac{1}{30};~~B_3=\frac{1}{42};~~
 B_4=\frac{1}{300},...
 \end{equation}

 The following lemma is needed:
  \begin{lemma}
 Let the relation: 
\begin{equation}\label{csum}
{\frac{1}{n!}=\sum_{k=0}^{n}c_k\frac{1}{(n-k+1)!}}
 \end{equation}
 Then  the constants $c_k$ of eq. (\ref{csum}) verifies the
 definition (\ref{defc}).
\end{lemma}

We need also another  formula similar to (\ref{csum}).
\begin{lemma}\label{constd}
Let the constants $d_k$ be defined by the relation: 
 
 \begin{equation}\label{dd}
 {\frac{1}{(n+2)!}=\sum_{k=0}^{n}d_k\frac{1}{(n-k+1)!}}~.
 \end{equation}

 Then the constants $c$ and $d$ are related by 
 
 \begin{equation}\label{cdd}
 {d_k= (-1)^kc_{k+1}}~.
 \end{equation}
 \end{lemma} 

Now we recall the main results established in \cite{sbcpol}:
\begin{Theorem}\label{bigtheorem}
Let $G$ be a semisimple Lie group.
If  $G_{\lambda}$ is  the generator
of the group $G$, then    ${\db G}_{\lambda}\in\D_1=\D_0\oplus
\D'_1$. More exactly, ${\db G}_{\lambda}\in {\A}_1$, i.e. 
  \begin{equation}\label{generic1}
{\db G}_{\lambda}= P_{\lambda}+\sum_{\beta \in\Delta_+} Q_{\lambda
 ,\beta} \pa_{\beta}, \lambda \in \Delta,
 \end{equation}
where $ P_{\lambda}$ and $ Q_{\lambda ,\beta}$  are polynomials in
 $z$. 

Explicitly:

 a) For  $\alpha\in \Delta_+$,

 \begin{equation}\label{plus}
 {\db  E}_{\alpha}=\sum_{k\ge 0}^{\nu}c_k\sum_{\beta\in\Delta_+}
 p_{k\alpha\beta}(z)\partial_{\alpha +\beta},
 \end{equation}
where  the coefficients  $c_k$, related to the Bernoulli numbers by eq.
 (\ref{cb}),
 are introduced  by eq. (\ref{csum}). The polynomials  $ p_{k\alpha\beta},
 k\in\N, \alpha\in\Delta_+$
  are given by the
 equation: 
 \begin{equation}\label{polp}
 p_{k\alpha\beta}(z)=
 \sum_{\stackrel{\alpha_1,\cdots ,\alpha_k}{\alpha_1+\cdots +\alpha_k =
 \beta}}n_{\alpha_1\cdots \alpha_k\alpha} z_{\alpha_1}\cdots z_{\alpha _k},
~k\ge 1
  \end{equation}
where
\begin{equation}\label{ndef}
 n_{\alpha_1\cdots \alpha_k\alpha}=
 n_{\alpha_1,\alpha}n_{\alpha_2,\alpha +\alpha_1}\cdots n_{\alpha_k,
 \alpha +\alpha_1+\cdots +\alpha_{k-1}},~(k\ge 1,\alpha_{0}=0), 
 \end{equation}
 and $n_{\alpha\beta}, \alpha ,\beta \in \Delta_+$ are the
 structure constants of eq. (\ref{cw}), and for $k=0$ the sum
(\ref{plus}) is just $\partial_{\alpha}$.  

 The expression (\ref{plus}) can be put also into a form in which the Bernoulli
 numbers are explicit:

\begin{equation}\label{plus1}
  {\db E}_{\alpha}=\partial_{\alpha}+
 \frac{1}{2}\sum_{\beta\in\Delta_+}z_{\beta}n_{\beta ,\alpha}
 \partial_{\alpha +\beta}+\sum_{k\ge 1}^{\nu}\frac{(-1)^{k-1}}{{(2k)!}}B_k
 \sum_{\beta\in \Delta_+}p_{k\alpha\beta}\partial_{\alpha +\beta}.
 \end{equation}

The degree of the polynomial $p$ has the property: $
 \mbox{\rm{degree~}} p_{k\alpha\beta}\le \nu;~
 p_{k\alpha\beta} ~$ as a function of $z$  contains only
 even powers. The table below contains the values of $\nu$.

\newpage
\begin{center}
 Degree $\nu$ for simple Lie algebras
 \end{center}
 $$\boxed{\begin{array}{llll} 
  A_l:\nu = l-1 & l\ge 1 & E_6: \nu = 10 & G_2: \nu =4\\
 B_l: \nu =2l-2 & l\ge 2  & E_7: \nu = 16 & \\
 C_l: \nu =2l-2 & l\ge 2  & E_8: \nu = 28    &\\
 D_l:\nu = 2l-4 & l\ge 3 & F_4: \nu = 10 & \\
\end{array}}$$ 

b) The differential action of the
  generators of the Cartan algebra is:

\begin{equation}\label{cartan}
 {\db H}= j +\sum_{\beta\in\Delta_+}\beta z_{\beta}\pa _{\beta}
 \end{equation}

c) If $(\alpha ,j)=0$, then

\begin{equation}\label{cartan1}
 {\db E}_{\alpha}=-\sum_{\beta\in\Delta_+}n_{\beta ,-\alpha}
 z_{\beta -\alpha} \partial_{\beta} .
 \end{equation}

d) If $\gamma \in\Delta_-$ is  simple root, then

 \begin{equation}\label{minus}
 {\db E}_{\gamma}= j \gamma z_{-\gamma} +\sum_{k\ge 0}^{\nu}d_k
 \sum_{\delta , \beta\in\Delta_+} q_{\gamma \delta}(z)
 p_{k\delta \beta}(z)\partial_{\beta +\delta},
 \end{equation}
where the coefficients $d$ are expressed through the coefficients $c$ by
 eq. (\ref{cdd}).

 The expression of the polynomials $q_{\gamma \delta},
 \gamma\in\Delta_- ,{\delta\in\Delta_+}$ is
 \begin{equation}\label{p1}
q_{\gamma \delta} = -\gamma z_{-\gamma}\delta
z_{\delta}  + \sum _{\mu\in\Delta_+}z_{\delta - \mu - \gamma} n_{\delta -\mu -
 \gamma , \gamma} z_{\mu}n_{\mu , \delta -\mu}.
 \end{equation}

 In the case of Hermitian symmetric cases
 eq. (\ref{plus}) becomes simply:
\begin{equation}
{\db E_{\alpha}}=\partial_{\alpha}
 \end{equation}
 while eq. (\ref{minus}) becomes
 \begin{eqnarray*}
 -{\db E}^-_{\alpha}& = &{\db K}^-_{\alpha}=(\alpha , j)z_{\alpha}+\frac{1}{2}
 z_{\alpha}\sum_{\beta\in\Delta_n^+}(\alpha ,\beta )z_{\beta}\partial_{\beta}
  \\
& & -\frac{1}{2}z_{\alpha}\sum_{\gamma -\alpha \in\Delta_k\setminus\{ 0\}}
 n_{\gamma ,-\alpha}
 n_{\gamma -\alpha ,-\beta}
 z_{\beta + \alpha -\gamma}\partial_{\beta} .
 \end{eqnarray*}
\end{Theorem}

The main ingredient 
 in the proof \cite{sbcpol} of  Theorem \ref{bigtheorem} is the formula:
\begin{equation}\label{bch1}
{\mbox{\rm{Ad}}}(\exp Z)=\exp {\mbox{\rm{ad}}}_Z,~ Z\in\got{g},
\end{equation}
i.e.
 \begin{equation}\label{bch}
 \e^ZX\e^{-Z}=\sum_{n \ge o} \frac{1}{n!}\ad^n_ZX,~ X, Z\in\got{g},
 \end{equation}
where
 $${\mbox{\rm ad}}_YX=[Y,X], ~~ 
 {\mbox{\rm ad}}^m_YX=[Y,{\mbox{\rm ad}}^{m-1}_Y X],~~ m>1, ~~\ad^0_YX=X.   $$ 
 We also have used the relation 
 \begin{equation}
 \e^Z\partial_{\alpha} (\e^{-Z})=-\left[\partial_{\alpha}(Z)+\sum_{n\ge 1}
 \frac{1}{(n+1)!}\ad^n_Z\partial_{\alpha}(Z)\right],
 \end{equation} 
 $$\partial_{\alpha}(Y)=\partial_{\alpha}Y-Y\partial_{\alpha}=
 -\ad_Y(\partial_{\alpha}) .$$
 
Due to the correspondence (\ref{correspond}),
 the operator of the left hand side of eq. (\ref{bch})
 corresponds to the
 differential action on \fl , and 
 $\e^Z\mb{X}\e^{-Z}\leadsto 
  -{\db X}~.$ 

\subsection{ An example: $SU(3)/\T$}\label{noul} 

In this section we follow closely \cite{saraceno} for the example
of the compact non-symmetric space  $M=SU(3)/\T$.

 The commutation relations of the generators are:
 
\begin{equation}\label{grun}
[C_{ij}, C_{kl}]=\delta_{jk}C_{il}-\delta_{il}C_{kj}, ~1\le i, j \le 3.
 \end{equation}

Let us consider the following parametrizations useful for the Gauss
decomposition and also in the definition of the coherent states for
the manifold $M$:
\begin{equation}\label{vzeta}
V_{+}(\zeta )=\exp
(\zeta_{12}C_{12}+\zeta_{13}C_{13}+\zeta_{23}C_{23}) ,
\end{equation}
\begin{equation}\label{vz}
V'_{+}(z )=\exp (z_{23}C_{23})\exp (z_{12}C_{12}+z_{13}C_{13}) .
\end{equation}

Let us denote by the same letter $C_{ij}$ the $n\times n$-matrix
having all elements $0$ except at the intersection of the line $i$
with the column $j$, that is $C_{ij}=(\delta_{ai}\delta_{bj})_{1\le a, 
b\le n}$; here  $n=3$. Then:
\begin{equation}\label{matricezeta}
V_{+}(\zeta )=
\left(
\begin{array}{ccc}
 1 & \zeta_{12} & \zeta_{13}+\frac{1}{2}\zeta_{12}\zeta_{23}\\
0 & 1 & \zeta_{23}\\
0 & 0 & 1
\end{array}
\right),
\end{equation}

\begin{equation}\label{matricez}
V'_{+}(z )=
\left(
\begin{array}{ccc}
 1 & z_{12} & z_{13}\\
0 & 1 & z_{23}\\
0 & 0 & 1
\end{array}
\right).
\end{equation}

Now observing that for 
\begin{equation}\label{egal}
z_{12}=\zeta_{12}; ~z_{13}=\zeta_{13}+\frac{1}{2}\zeta_{12}\zeta_{23};~
z_{23}=\zeta_{23} , 
\end{equation}
we get 
\begin{equation}\label{eegal}
V_{+}(\zeta )=V'_{+}(z ) .
\end{equation}
So, we have two  parametrizations  of the compact non-symmetric flag manifold
 $M=SU(3)/S(U(1))\times U(1)\times U(1))$: one in  $\zeta$, given by eq. 
(\ref{matricezeta}) and the other one in $z$, given by  (\ref{matricez}), which
 are identified with the relations (\ref{egal}).

Let us consider also the CS-vectors
\begin{equation}\label{vect}
\begin{array}{ccl}
\phi_z & = & [V'_{+}(z )]^+\phi_{w}\\
& = & \exp (\bar{z}_{12}{\mb{C}}_{21}+\bar{z}_{13}\mb{C}_{31})
\exp (\bar{z}_{23}\mb{C}_{32})\phi_{w}.
\end{array}
\end{equation}

Now $\phi_w$ is chosen as maximal weight vector corresponding to the weight
$w=(w_1,w_2,w_3)$
 such that  $ j_1=\omega_1-\omega_2\ge 0, j_2=
 \omega_2-\omega_3\ge 0 $ and the lowering operators are $\mb{C}_{ij}, i>j$, 
and $C_{ii}$ generates  the Cartan algebra, i.e.
\begin{equation}\label{altele}
\left\{
\begin{array}{cll}
{\mb C}_{ij}\phi_w & \not=&  0,~ i>j ,\\
{\mb C}_{ij}\phi_w & = & 0 ,~i<j ,\\
{\mb C}_{ii}\phi_w & = & w_i\phi_w .
\end{array}
\right.
\end{equation}
The coherent state vectors corresponding to the representation
$\pi_w$ determined by eqs. (\ref{altele})  are introduced as
\begin{equation}\label{noucoer}
e_z=\pi_{w}((V'_{+}(z ))\phi_w .
\end{equation}
Denoting by $Z$ the matrix
 \begin{equation}\label{zmat} 
Z=
 \left(
 \begin{array}{lll}
 1 & z_{12} & z_{13} \\
 0 & 1 & z_{23} \\
 0 & 0 & 1
 \end{array}
 \right) ,
\end{equation}
the reproducing kernel which determines the scalar product
$(e_{\bar{z}},e_{\bar{z}})$
has the expression: 
\begin{eqnarray*}
 K(ZZ^+) & = &\Delta^{j_1}_1(ZZ^+)\Delta^{j_2}_2(ZZ^+);\\
 \Delta_1 & = & 1+|z_{12}|^2+|z_{13}|^2;\\
\Delta_2 & = & (1+|z_{12}|^2+|z_{13}|^2)(1+|z_{23}|^2)
-|z_{12}+z_{13} \overline{z}_{23}|^2 .
 \end{eqnarray*}

In particular, it is observed that for 
$z_{23}= 0$, $M$ becomes $ SU(3)/S(U(2)\times U(1))=G_1(\C^3)=\C\P^2$.
In order to compare with the scalar product for  coherent states on
$M\approx\C\P^2\approx G_1(\C^3)$ we remember that in the case of the
Grassmannian we have used in \cite{sbcag,sbl}
 a weight which here corresponds to $w_1=1,
w_2=w_3=0$
and then on $\C\P^2$ the reproducing kernel is just
$K(ZZ^+)= \Delta_1$.

The calculation which has as result  Lemma \ref{ceurile}
(cf. \cite{sbcpol}) does 
not use the value of the reproducing kernel and is an algebraic one.

 \begin{lemma}\label{ceurile}
The differential operators ${\db{C}}_{ij}$
associated to  the generators $C_{ij}$ are given by the formulas:
  \begin{eqnarray*}
{\db C}_{11} & = & -z_{12}\pa_{12}-z_{13}\pa_{13}+ w_1 ,\\
  {\db C}_{12} & = &  \pa_{12} ,\\
 {\db C}_{13} & = & \pa_{13} ,\\
 {\db C}_{21} & = & -z_{12}^2\pa_{12}-z_{12}z_{13}\pa_{13}+
 (z_{12}z_{23}-z_{13})\pa_{23}+(w_1-w_2)z_{12} ,\\ 
 {\db C}_{22} & = & z_{12}\pa_{12}-z_{23}\pa_{23}+w_2 ,\\
 {\db C}_{23} & = & z_{12}\pa_{13}+\pa_{23} ,\\
 {\db C}_{31} & = &
 -z_{12}z_{13}\pa_{12}-z^2_{13}\pa_{13}+
 \underline{(z_{12}z_{23}-z_{13})z_{23}}\pa _{23}+\\
 & & (w_1-w_3)z_{13}-(w_2-w_3)z_{12}z_{23} ,\\
 {\db C}_{32} & = & z_{13}\pa_{12}-z_{23}^2\pa_{23}+(w_2-w_3)z_{23} ,\\
 {\db C}_{33} & = & z_{13}\pa_{13}+z_{23}\pa_{23} + w_3~.
  \end{eqnarray*}
\end{lemma}

 We have underlined the apparition of a third-degree polynomial multiplying the
 partial derivative of ${\db C}_{31}$.
 Note also the relation
$C_{11}+C_{22}+C_{33}=w_{1}+w_{2}+w_{3}$\\[3ex]

\section{ Equations of motion}\label{eqm}

Let $(M\approx G/H,\omega )$ be a
 quantizable homogeneous K\"ahler manifold. Let us consider the triplet
 $(L,h,\nabla )$, where $L$ is a $G$-homogeneous
 holomorphic positive line bundle, 
$h$ the  hermitian metric and   $\nabla$ the  connection compatible with 
the complex structure.
 The quantization  condition reads: 
\begin{equation}\label{quant}
\omega=i\Theta_L=\pi c_1(L)=-i\pa\bar{\pa} \log (h), 
\end{equation}
where  $\Theta_L$ is  the curvature matrix of  $L$ and 
$c_1$ is   first Chern class. Here the hermitian metric is given by  
\begin{equation}\label{hh}
h(z,z)=(e_{\bar{z}},\e_{\bar{z}})^{-1},
\end{equation}
 while the  K\"ahler potential  is 
$K(z)=-\log h(z)= \log(\e_{\bar{z}},\e_{\bar{z}})$. The K\"ahler
two-form is 
\begin{equation}\label{twoform}
\omega (z)
 = i\sum_{\alpha ,\beta\in\Delta_+}\NC_{\alpha,\beta}dz_{\alpha} \wedge
d\bar{z}_{\beta},
\end{equation}
where
\begin{equation}\label{pot}
\NC_{\alpha,\beta}(z)=\frac{\pa^2}{\pa z_{\alpha}\pa\bar{z}_{\beta}}\log
(e_{\bar{z}},e_{\bar{z}}). 
\end{equation}

If $H\in\got{g}$ and the corresponding element in 
$\am $ is $\mb{H}$, then the  energy function (covariant symbol) 
attached to  $\mb{H}$ is 
 \begin{equation}\label{energie}
\HR (z,\bar{z}):=(e_{\sigma (z)},\mb{H}e_{\sigma (z)})
=N^{2}(z)(e_{\bar{z}},\mb{H}e_{\bar{z}})
=\frac{(e_{\bar{z}},\mb{H}e_{\bar{z}})}{(e_{\bar{z}},e_{\bar{z}})} .
\end{equation}

The homogeneous quantizable K\"ahler manifolds $M\approx G/H$ for which
the group $G$ verify some  obstructions which for semisimple
Lie groups are
$H^1(\got{g})=H^2(\got{g})=\{0\}$
can be organized as elementary Hamiltonian
$G$-spaces \cite{gs}.
  Passing on from the dynamical system problem
 in the Hilbert space $\Hi$ to the corresponding one on $M$ is called
sometimes {\it dequantization}, and the system on $M$ is a classical
one. Following Berezin \cite{berezin1}, the
motion on the classical phase space can be described by the local
equation of motion 
\begin{equation}\label{ecmiscare}
\dot{z}_{\gamma}=i\{\HR ,z_{\gamma}\},\gamma\in\Delta_+ .
\end{equation}
In eq. (\ref{ecmiscare}) $\{\cdot,\cdot \}$ denotes  the  Poisson bracket: 
$$\{ f,g\} =\sum_{\alpha ,\beta\in\Delta_+}\NC^{-1}_{\alpha
,\beta}\left\{\frac{\pa f}{\pa
z_{\alpha}}\frac{\pa g}{\pa\bar{z}_{\beta}}-\frac{\pa f}{\pa\bar{z}_{\alpha}}
\frac{\pa g}{\pa z_{\beta}}\right\}, f,g\in C^{\infty}(M) .$$
The equations of motion (\ref{ecmiscare}) can be written down as
\begin{equation}
i\left(
\begin{array}{cc}
0 & \NC\\
 & \\
-\bar{\NC} & 0
\end{array}\right)
\left(
\begin{array}{c}
\dot{z} \\
  \\
\dot{\bar{z}}
\end{array}
\right)
=-\left(
\begin{array}{c}
\frac{\pa}{\pa z}\\
   \\
\frac{\pa}{\pa \bar{z}}
\end{array}
\right)\HR .
\end{equation}

We consider algebraic Hamiltonian linear in the generators of the group
\begin{equation}\label{hsum}
\mb{H}=\sum_{\lambda\in\Delta}\epsilon_{\lambda}{\mb{G}}_{\lambda} .
\end{equation}
If we take into account eq. (\ref{generic}),
then the {\em  equations of motion (\ref{ecmiscare}) are} (cf. \cite{sbcag,sbl})
\begin{equation}\label{move}
{i\dot{z}_{\alpha}=\sum_{\lambda\in\Delta}\epsilon_{\lambda}Q_{\lambda 
,\alpha}, \alpha\in\Delta_+.}
\end{equation}

For linear Hamiltonians (\ref{hsum}), we look for the
solution of the Schr\"odinger equation (\ref{sch})
\begin{equation}\label{sch}
\mb{H}\psi=i\dot{\psi}
\end{equation}
 via Perelomov's coherent state vectors:
\begin{equation}\label{coerentesr}
\psi=\psi (z)=e^{i\varphi}e_z .
\end{equation}
 It can be proved (cf. \cite{sbcag,sbl}) that:  {\it If the initial
state is a coherent state and the Hamiltonian is linear in the
generators of the group, then the state will evolve into a coherent
state}. More exactly, 
\begin{Proposition}\label{ecmisc}
On the manifold $M$ of coherent states for which formulas
(\ref{generic}) are true, the classical motion and the quantum evolution 
generated by the linear Hamiltonian (\ref{hsum}) are given both
by the same 
equation of motion (\ref{move}). For semisimple Lie groups the
expression
of the polynomials $Q$-s is that given in Theorem \ref{bigtheorem}.
 The phase in eq. (\ref{coerentesr})
is given by the sum (\ref{unghitotal1})
\begin{equation}\label{unghitotal1}
\varphi = -\int_{0}^{t}\HR dt -\Im\int_{0}^{t}
\frac{(e_{\bar{z}},de_{\bar{z}})}{(e_{\bar{z}},e_{\bar{z}})} 
\end{equation}
 of the dynamical and Berry phase.
\end{Proposition}

\section{Explicit examples of equations of motions}

The first two examples presented below are pedagogical. The case of
the Grassamnn manifold and its non-compact dual
are taken from \cite{sbcag,sbl}. The last
example is derived as application of the expressions of the
differential
operators presented in \S \ref{noul}.

 \subsection{ The oscillator Group}\label{osc} 
The canonical commutation relations of the creation and annihilation
operators are
$$[\mb{a}_{\lambda},\mb{a}^+_{\lambda '}]=\delta_{\lambda ,\lambda '}1; \lambda
,\lambda '=1,\dots , n .$$
The Perelomov's coherent state vectors (Glauber's coherent states) are
$$e_z=e^{\sum z_{\lambda}\mb{a}^+_{\lambda}}e_0 , $$
where 
$$\mb{a}_{\lambda}e_0=0 ,$$
$$\mb{a}^+_{\lambda}e_z=\pa_{\lambda}e_z ,$$
$$\mb{a}_{\lambda}e_z=z_{\lambda}e_z .$$
The differential operators are 
${\db G}_{\lambda ,\mu}=P_{\lambda ,\mu}+\sum Q_{\lambda ,\mu
;\beta}\pa _{\beta} ,$
where 
$P_{\lambda ,\mu}=0;~ Q_{\lambda ,\mu ;\beta}=z_{\mu}\delta_{\beta
,\lambda} .$

The linear Hamiltonian  
$$\mb{H}=\sum \omega_{\lambda
,\mu}\mb{a}^+_{\lambda}\mb{a}_{\mu}+f_{\lambda}\mb{a}^+_{\lambda}
+ \bar{f}_{\lambda}a_{\lambda}$$
implies the  equation of motion
$$
{i\dot{z}_{\alpha}=\omega_{\alpha ,\mu}z_{\mu}+f_{\alpha}} .$$

\subsection{ $SU(2)/U(1)$ \& $SU(1,1)/U(1)$}\label{sfera} 
The generators verify the commutation relations
$$[J_0,J_{\pm}]=\pm J_{\pm}; [J_-,J_+]=-2J_0 ,$$
and 
$$\left\{
\begin{array}{ccl}
\mb{J}_+ e_{j,-j} & \not= & 0 ,\\
\mb{J}_- e_{j,-j} & = & 0 ,\\
\mb{J}_0e_{j,-j} & = & -je_{j,-j} .
\end{array}
\right.$$
We have
$$\left\{
\begin{array}{ccl}
\mb{J}_0e_z & = & (-j+z\pa )e_z ,\\
\mb{J}_+e_z & = & \pa e_z ,\\
\mb{J}_-e_z & = & (2jz-z^2\pa )e_z .
\end{array}
\right.$$
The linear Hamiltonian
$$\mb{H}=\epsilon_+\mb{J}_++\epsilon_-\mb{J}_-+\epsilon_0\mb{J}_0 ,$$
where
$\epsilon_0^+=\epsilon_0 ;~
\epsilon^+_+=\epsilon_-;~\epsilon^+_-=\epsilon_+ ,$
implies the Riccati equation of motion
\begin{equation}\label{aro}
{i\dot{z}=\epsilon_0z-\epsilon_-z^2+\epsilon_+}.
\end{equation}

In the non-compact case, the generators verifies the commutation relations
$$[K_0,K_{\pm}]=\pm K_{\pm}; [K_-,K_+]=2K_0$$
and finally we have the equation of motion with a minus sign in the
front of the term $z^2$ comparatively with the equation (\ref{aro})
corresponding to the compact case:
\begin{equation}\label{aro1}
{i\dot{z}=\epsilon_0z+\epsilon_-z^2+\epsilon_+}.
\end{equation}
The change of sign for the Riemann sphere and his non-compact dual as
in equations (\ref{aro}) and (\ref{aro1}) take place also in the case
of the Grassmann manifold.

\subsection{ The complex Grassmannian \Gras}\label{grasu} 
The complex Grassmann manifold is
$\Gras\approx SU(n+m)/S((U(n)\times U(m)) . $
The non-compact dual of the complex Grassmann manifold is
$SU(n,m)/S((U(n)\times U(m))$.

The differential operators 
on  \Gras~ (``Slater determinant
manifold'') are (cf. \cite{sbcag,sbl}):
$$\left\{
\begin{array}{l}
\db{E}^+_{im}  =  {\db K}^+_{im} =  \pa_{im} ,\\
-\db{E}^-_{mi}  =  \db{K}^-_{mi}
=(j_i-j_m)z_{im}+\sum z_{jm}z_{in}\pa_{jn} ,\\
\db{H}_{\mu\nu} = \delta_{\mu\nu}j_{\nu}-z_{im}
(\delta_{m\mu}\pa_{i\nu}-\delta_{i\nu}\pa_{m\mu}) . 
\end{array}
\right.
$$
The linear Hamiltonian is
\begin{eqnarray*}
 H & = & \sum_{1\le\mu ,\nu \le m}(\ep^0_1 )_{\mu\nu}H_{\mu\nu}
+\sum_{m<\mu ,\nu\le m+n}(\ep ^0_2)_{\mu\nu}H_{\mu\nu}\\
& + &
\sum_{i=1}^m\sum_{p=m+1}^{m+n}\ep^+_{ip}F^+_{ip}
+\ep^-_{ip}F^-_{pi} ,
\end{eqnarray*}
where $F$ $(K)$ are   generators for compact (resp. non-compact)
Grassmannian. \\
{\it The equation of motion on the compact (resp., non-compact)
Grassmann manifold is   the Matrix Riccati equation} (cf. \cite{sbcag,sbl})
\begin{equation}\label{Rmatr}
{ i\dot{Z}=-Z\ep ^0_2 + \ep^0_1Z +\ep^+\pm Z\ep^-Z} , 
\end{equation}
$$(\ep^0_{1,2})^+=\ep^0_{1,2}; (\ep^+)^+=\ep^- .$$
With the linear fractional change of variables
$Z=XY^{-1}$ we get a linearization
of  the matrix Riccati equation (\ref{Rmatr})\\
\begin{equation}\label{ric}
\left(
\begin{array}{c}
\dot{X}\\
\dot{Y}
\end{array}\right)
=h_{n,c}
\left(
\begin{array}{c}
X\\
Y
\end{array}\right),
~h_{n,c}=
\left(
\begin{array}{cc}
-i \ep^0_1 & -i\ep^+ \\
\pm\ep^- & i\ep ^0_2
\end{array}\right) ,
\end{equation}
where the subindex $c$ ($n$) corresponds to the  compact (resp.,
  non-compact) case. 
The dimension of the matrices $ X, Y, Z, \ep^0_1, \ep ^0_2, \ep^+,
\ep^-$ are, respectively:  $m\times n, n\times n, m\times n, m\times m, 
n\times n, m\times n$, and $n\times m$.

\subsection{  $SU(3)/\T$}\label{ultim} 

Now we come back to the example of \S \ref{noul} and 
 we consider the
 linear hermitian Hamiltonian:

\begin{equation}\label{ham3}
{\mb H}= \sum_{i,j=1}^3\ep_{ij}{\mb C}_{ij};  ~~\ep_{ij}^+=\ep_{ij}.
\end{equation}

\begin{Proposition}\label{misc3}
In the parametrization (\ref{matricez})
on  $SU(3)/S(U(1)\times U(1)\times U(1))$,
 the equations of motion (\ref{move})
associated to the Hamiltonian (\ref{ham3}) are:
\begin{equation}\label{treiec}
{
\begin{array}{ccl}
i\dot{z}_{12}&=&
-\ep_{11}z_{12}+\ep_{12}-\ep_{21}z_{12}^2+\ep_{22}z_{12}-\ep_{31}z_{12}z_{13}
+\ep_{32}z_{13} \\
i\dot{z}_{13}&=&-\ep_{11}z_{13}+\ep_{13}-\ep_{21}z_{12}z_{13}+
\ep_{23}z_{12}-\ep_{31}z_{13}^2+\ep_{33}z_{13} \\
i\dot{z}_{23}&=& \ep_{21}(z_{12}z_{23}-z_{13})-\ep_{22}z_{23}+\ep_{23} 
+\ep_{31}(z_{12}z_{23}-z_{13})z_{23}\\
& & -\ep_{32}z_{23}^2+\ep_{33}z_{23}\\
\end{array}
}
\end{equation}
\end{Proposition}

{\bf Comment}
In the 
$z$-parametrization given by eq. (\ref{matricez})
 the first two equations in (\ref{treiec})  do not depend on $z_{23}$
and are in fact a  matrix Riccati equation on 
$ SU(3)/S(U(2)\times U(1))=G_1(\C^3)=\C\P^2$, in accord with
eq. (\ref{Rmatr}).  In the situation
$z_{23}=0$, the representations in $z$ and $\zeta$ coincides, and  
the terms in  Proposition (\ref{misc3})
corresponding to $\frac{\pa}{\pa z_{12}} $
are missing and we get in the equation of motion (\ref{treiec})
  only the first two
equations.
 In the parametrization with $z$
  the flag structure of $SU(3)/S(U(1))\times U(1)\times U(1))$  is
reflected in the decoupling of the first two equations of motion
of the third.
\footnotesize

{\bf Acknowledgment} S. B. express his thanks to the Organizers of
the Fifth International Workshop of Differential Geometry and its
Applications for the opportunity to present this talk in
Timisoara. Discussions during the Workshop with Professor Lieven
Vanhecke are kindly acknowledged. 
S.B. is grateful to Professor Karl-Hermann Neeb  for many
suggestions and criticism. 

\end{document}